\documentclass{birkjour}
\usepackage{amsfonts,amscd}
\usepackage{amssymb}
\usepackage{url}
\usepackage[english]{babel}
\theoremstyle{plain}
\newtheorem{theorem}                {Theorem}      [section]

\newtheorem{lemma}        [theorem]  {Lemma}

\theoremstyle{definition}
\newtheorem{example}      [theorem]  {Example}
\newtheorem{remark}       [theorem]  {Remark}

\newcommand{\secref}[1]{\S\ref{#1}}
\numberwithin{equation}{section}

\def \R{{\mathbb R}}

\def \z{{\mathbb Z}}

\usepackage{color}

\numberwithin{equation}{section}

\begin{document}

\title[Equivariant biharmonic maps]{A general approach to equivariant biharmonic maps}

\author{S.~Montaldo}
\address{Universit\`a degli Studi di Cagliari\\
Dipartimento di Matematica e Informatica\\
Via Ospedale 72\\
09124 Cagliari, Italia}
\email{montaldo@unica.it}

\author{A.~Ratto}
\address{Universit\`a degli Studi di Cagliari\\
Dipartimento di Matematica e Informatica\\
Viale Merello 93\\
09123 Cagliari, Italia}
\email{rattoa@unica.it}

\begin{abstract}
In this paper we describe a 1-dimensional variational approach to the analytical construction of
 equivariant biharmonic maps. Our goal is to provide a direct method which enables analysts to compute directly the analytical conditions
 which guarantee biharmonicity in the presence of suitable symmetries. In the second part of our work, we illustrate and discuss some examples. In particular, we obtain a 1-dimensional stability result, and also show that biharmonic maps do not satisfy the classical maximum principle proved by Sampson for harmonic maps.
\end{abstract}




\thanks{Work supported by Contributo d'Ateneo, University of Cagliari, Italy.}

\maketitle

\section{Introduction}\label{intro}
{\it Harmonic maps}  are critical points of the {\em energy} functional
\begin{equation}\label{energia}
E(\varphi)=\frac{1}{2}\int_{M}\,|d\varphi|^2\,dv_g \,\, ,
\end{equation}
where $\varphi:(M,g)\to(N,h)$ is a smooth map between two Riemannian
manifolds $M$ and $N$. In analytical terms, the condition of harmonicity is equivalent to the fact that the map $\varphi$ is a solution of the Euler-Lagrange equation associated to the energy functional \eqref{energia}, i.e.
\begin{equation}\label{harmonicityequation}
    {\rm trace} \, \nabla d \varphi =0 \,\, .
\end{equation}
The left member of \eqref{harmonicityequation} is a vector field along the map $\varphi$, or, equivalently, a section of the pull-back bundle $\varphi^{-1} \, (TN)$: it is called {\em tension field} and denoted $\tau (\varphi)$. In local charts, the tension field, which is the trace of the second fundamental form, is given by the following expression:
\begin{equation}\label{tensionfield}
    \tau^\gamma (\varphi)= g^{ij} \left ( \nabla(d \varphi) \right )_{ij}^\gamma \,\, ,
\end{equation}
where
\begin{equation}\label{secondaformalocale}
    \left ( \nabla(d \varphi) \right )_{ij}^\gamma = \frac{\partial^2 \varphi^\gamma}{\partial x^i \, \partial x^j} - {}^M\Gamma_{ij}^k \, \frac{\partial \varphi^\gamma}{\partial x^k} \, + \, {}^N\Gamma_{\alpha \beta}^\gamma \, \frac{\partial \varphi^\alpha}{\partial x^i} \, \frac{\partial \varphi^\beta}{\partial x^j} \,\, .
\end{equation}
In \eqref{secondaformalocale}, ${}^M\Gamma$ and ${}^N\Gamma$ denote respectively the Christoffel symbols of the Levi-Civita connections of $(M,g)$ and $(N,h)$; note also that Einstein's convention of sum over repeated indices is adopted. In particular, inspection of \eqref{secondaformalocale} reveals that the harmonicity equation \eqref{harmonicityequation} is a second order, elliptic system of partial differential equations.

Every vector field $v$ along $\varphi$ defines a variation of $\varphi$ by setting
\begin{equation}\label{variation}
    \varphi_t (p)= \exp_{\varphi(p)}\left( t\,v(p)\right) \,\, , \qquad t\in(-\varepsilon ,\varepsilon)
\end{equation}
(note that $\varphi_0=\varphi$). If $v$ is compactly supported, then
\begin{equation}\label{decrescenzarapidaenergia}
  \nabla_v \, E(\varphi)= \left . \frac{d\,E(\varphi_t) }{dt}\right |_{t=0} =-\, \int_{M}\,\langle\tau(\varphi) ,\, v\rangle\,dv_g \,\, .
\end{equation}
In particular, it follows that \eqref{harmonicityequation} is equivalent to the vanishing of the directional derivative in \eqref{decrescenzarapidaenergia} for all $v$. Moreover, we point out that \eqref{decrescenzarapidaenergia} implies the fact that the tension field $\tau (\varphi)$ provides the direction in which the energy decreases more rapidly.

The study of harmonic maps is a very wide area of research, involving a rich interplay of geometry,
analysis and topology. We refer to \cite{PBJCW,JEELLE2,JEELLE1,Xin} for notation and
background on harmonic maps and to \cite{BMBib} for a more recent
bibliography.

A related topic of growing interest deals with the study of the so-called {\it biharmonic maps}: these maps, which provide a natural generalisation of harmonic maps, are the critical points of the bienergy functional (as suggested by Eells--Lemaire \cite{EL83})
\begin{equation}\label{bienergia}
    E_2(\varphi)=\frac{1}{2}\int_{M}\,|\tau (\varphi)|^2\,dv_g \,\, .
\end{equation}
In \cite{Jiang} Jiang derived the first variation and the second variation formulas for the bienergy. In particular, he showed that the Euler-Lagrange equation associated to $E_2(\varphi)$ is
\begin{equation}\label{bitensionfield}
    \tau_2(\varphi) = - J\left (\tau(\varphi) \right ) = - \triangle \tau(\varphi)- \rm{trace} R^N(d \varphi, \tau(\varphi)) d \varphi = 0 \,\, ,
    \end{equation}
where $J$ is (formally) the Jacobi operator of $\varphi$, $\triangle$ is the rough Laplacian defined on sections of $\varphi^{-1} \, (TN)$ and
\begin{equation}\label{curvatura}
    R^N (X,Y)= \nabla_X \nabla_Y - \nabla_Y \nabla_X -\nabla_{[X,Y]}
\end{equation}
is the curvature operator on $(N,h)$.

In this context, the biharmonic version of \eqref{decrescenzarapidaenergia} is (see\cite{Jiang})
\begin{equation}\label{decrescenzarapidabienergia}
  \nabla_v \, E_2(\varphi)= \left . \frac{d\,E_2(\varphi_t) }{dt}\right |_{t=0} = \int_{M}\,\langle  \tau_2(\varphi) ,\, v\rangle \,dv_g \,\, .
\end{equation}
Therefore, \eqref{bitensionfield} represents the vanishing of the directional derivative in \eqref{decrescenzarapidabienergia} for all $v$. In particular, the bitension field $\tau_2 (\varphi)$ provides the direction in which the bienergy decreases more rapidly.

We point out that \eqref{bitensionfield} is a {\it fourth order}  semi-linear elliptic system of differential equations. We also note that any harmonic map is an absolute minimum of the bienergy, and so it is trivially biharmonic; thus, a general working plan is to study the existence of biharmonic maps which are not harmonic that we shall call {\it proper biharmonic}. We refer to \cite{Montaldo} for existence results and general properties of biharmonic maps.

Fourth order differential equations are of great importance in various fields. By way of example, we cite an instance which is well-known to civil engineers: the structural problem of a beam resting on an elastic soil is amenable to the following differential equation, in which the unknown function $Y(x,t)$ represents the response of the beam in position $x$, at the time $t$:
\begin{equation}\label{trave}
    \frac{\partial^2}{\partial x^2} \left( EI(x) \, \frac{\partial^2 \, Y }{\partial x^2}\right) + m \, \frac{\partial^2 \, Y}{\partial t^2} + k \, Y= f(x,t) \,\, ,
\end{equation}
where $EI(x)$ measures the flexural rigidity (Young modulus), $m$ and $k$ are two constants depending respectively on the material of the beam and on the elasticity of the soil; the right member $f(x,t)$ represents the external load (see \cite{Corradi} for details).

In general, we encounter enormous difficulties to study fourth order differential equations. In particular, the presently known instances of biharmonic maps have always been obtained essentially by means of geometric intuition and simplification (see, for example, \cite{BMO1,BMO2}).

In this order of ideas, {\it equivariant theory} deals with special families of maps having enough symmetries  to guarantee that harmonicity reduces to the study of a second order {\it ordinary}
differential equation (we refer to \cite{EELRAT, Xin} for background and examples).

In this paper we shall explain how this general framework can be adapted to include the study of biharmonic maps. Although the study of the resulting fourth order ordinary differential equation remains, in general, a very difficult task, one of our aims is to provide analysts with a direct approach which permits them the computation of the relevant equations avoiding, in particular, the often heavy burden of dealing with an advanced riemannian geometric machinery. As an application, we obtain a 1-dimensional stability result and also show that proper biharmonic maps in general do not satisfy Sampson's maximum principle.

\noindent{\bf Acknowledgement}. 
The authors wish to thank the referee for comments 
and suggestions that have improved the paper.

\section{The 1-dimensional variational approach}\label{sezvarapproach}

In order to introduce the variational context in which we work, let us first recall some basic facts concerning equivariant maps and their associated
1-dimensional variational problem. Various, rather general approaches are possible (see, for example, \cite{PBJCW, EELRAT,Xin}): here we adopt the setting of \cite{Xin}, which is largely sufficient for our purposes and, at the same time, will ease our task to make this article as self-contained as possible. In particular, we shall consider maps $f$ which are \emph{equivariant with respect to Riemannian submersions}: keeping, just for this short introduction to equivariance, the notation of \cite{Xin}, that amounts to require that the following be a commutative diagram:
\begin{equation}\label{diagrammacommutativo}
\begin{CD}
E_1 @>f>>E_2\\
@V \pi_1VV @V V\pi_2V\\
M_1 @>\bar{f}>> M_2.
\end{CD}
\end{equation}
where $\pi_i : E_i \to M_i \,, \,\, i=1,2$ are Riemannian submersions.
We say that a vector field $v$ along $f$ is \emph{basic} if: for any $p\in E_1$, $v(p)$ is horizontal with respect to $\pi_2$ and there exists a vector field $\bar{v}$ along $\bar{f}$ such that
\begin{equation}\label{basicvectorfield}
d\pi_2(v(p))=\bar{v}(\pi_1(p)) \,\, .
\end{equation}
We need the following adaptation of a result of Xin (see Theorem 6.5 in \cite{Xin}) to the biharmonic case:
\begin{lemma}\label{adaptationofXin}
Let $f:E_1\to E_2$ be an equivariant map with respect to Riemannian submersions. If its bitension field
$\tau_2$ is a basic vector field, then $f$ is biharmonic if and only if it is a critical point of the bienergy with respect to all equivariant variations.
\end{lemma}
\begin{proof} Although the proof follows \cite{Xin} step by step, we report it here for the sake of completeness. We only have to prove that if $f$ is a critical point with respect to all equivariant variations then $f$ is biharmonic.
Let $h$ be any equivariant function on $E_2$ with compact support and $\bar{h}$ its induced function on $M_2$. Consider the following variational maps
$$
f_t(p)=\exp _{f(p)}(t\, h(f(p))\, \tau_2(f))\,\,,\quad t\in(-\varepsilon,\varepsilon).
$$
By construction the corresponding variational vector field is $(h\circ f)\tau_2(f)$. We now show that $f_t$ is an equivariant variation of $f$. Since $\tau_2(f)$ is basic, there exists a vector field $\bar{v}(\bar{f})$ along $\bar{f}$  such that
$d\pi_2(\tau_2(f)(p))=\bar{v}(\bar{f})(\pi_1(p))$. Thus, we have
\begin{eqnarray*}
\pi_2\circ f_t(p)&=&\exp _{\pi_2\circ f(p)}(d\pi_2(t\, h(f(p))\, \tau_2(f)))\\
&=& \exp _{\bar{f}(\pi_1(p))}(t\, \bar{h}\circ \bar{f}(\pi_1(p))\, \bar{v}(\bar{f})(\pi_1(p)))\\
&=& \bar{f}_t(\pi_1(p))\,,
\end{eqnarray*}
with
$$
\bar{f}_t(q)=\exp_{\bar{f}(q)}(t\, \bar{h}\circ \bar{f} (q)\, \bar{v}(\bar{f})(q))\,,\quad q\in M_1\,,\, t\in(-\varepsilon,\varepsilon)\,,
$$
where $\bar{h}$ is such that $h(p)=\bar{h}(\pi_2(y))$, for any $y\in E_2$.
Now, by using the first variational formula \eqref{decrescenzarapidabienergia}, we obtain
\begin{eqnarray*}
0=\left . \frac{d\,E_2(f_t) }{dt}\right |_{t=0} &=& \int_{E_1}\,\langle \tau_2(f) ,\, (h\circ f)\tau_2(f)\rangle \,dv_g \\
&=& \int_{E_1}\,(h\circ f) \, |\tau_2(f)|^2 \,dv_g\,,
\end{eqnarray*}
which implies, by the assumptions, that $\tau_2(f)\equiv 0$.
\end{proof}

In order to describe our applications, we restrict our attention to the case that the base manifolds in \eqref{diagrammacommutativo} are 1-dimensional. More specifically, we consider equivariant maps $\varphi_{\alpha}:(M,g)\to(N,h)$ and denote by
\begin{equation}\label{alfa}
\alpha : [a,b] \to \R
\end{equation}
the associated map between the base manifolds. Note that this function $\alpha$ may have to satisfy suitable boundary conditions dictated by the geometry of the problem.  If these maps have enough symmetries (in particular, to ensure that the tension field is basic) (see \cite{EELRAT} for details, and \secref{sezesempi} for examples), then the energy functional \eqref{energia} takes the following form:
\begin{equation}\label{energiaridotta}
    E^{\varphi}(\alpha)= \int_a^b \,\, L(t, \alpha, \dot{\alpha}) \, dt \,\, ,
\end{equation}
where $\dot{\alpha}$ denotes derivative with respect to $t$, and $L(\cdot, \cdot, \cdot)$ is a differentiable function depending on the geometry of the problem under consideration. Now, according to (the harmonic maps version of) Lemma \ref{adaptationofXin}, the harmoni\-city of the map $\varphi_{\alpha}$ is equivalent to the fact that $\alpha $ is a {\it critical point} of the so-called reduced energy functional \eqref{energiaridotta}. That is to say, the function $\alpha$ must be a solution of the Euler-Lagrange equation associated with \eqref{energiaridotta}, which is the following:
\begin{equation}\label{eqeulero1}
    \frac{\partial L}{\partial \alpha} - \frac{d}{dt} \,  \left (\frac{\partial L}{\partial \dot{\alpha}} \right ) \, = \, 0 \,\, .
\end{equation}
For future reference, it is useful to recall how equation \eqref{eqeulero1} is derived from \eqref{energiaridotta}. The first step in this direction is to recognize that the requirement that $\alpha$ be critical is equivalent to the vanishing of the directional derivative
\begin{equation}\label{derivdir}
    \frac{d}{dh} \, \left [  E^{\varphi}(\alpha \, +\, h\, \beta) \right ] |_{h=0}
    \end{equation}
for all compactly supported (differentiable) variations $\beta : [a,b] \to \R$. Next, we compute \eqref{derivdir} explicitly:
\begin{eqnarray} \label{calcolo1}
\nonumber
  \frac{d}{dh} \, \left [  E^{\varphi}(\alpha \, +\, h\, \beta) \right ] |_{h=0} &=& \frac{d}{dh}\left ( \int_a^b \,\, L(t, \alpha+h\beta, \dot{\alpha}+h\dot{\beta})\, dt  \right )\Big {|}_{h=0} \\ \nonumber
   &=& \int_a^b \,\left (\frac{d}{dh}\, L(t, \alpha+h\beta, \dot{\alpha}+h\dot{\beta})|_{h=0} \right )\, dt  \\ \nonumber
   &=& \int_a^b \,\left (  \frac{\partial L}{\partial \alpha} \, \beta +  \frac{\partial L}{\partial \dot{ \alpha}} \, \dot{\beta}\right ) \, dt \\
   &=& \int_a^b \left [ \left ( \frac{\partial L}{\partial \alpha} - \frac{d}{dt} \,  \frac{\partial L}{\partial \dot{\alpha}} \right ) \, \beta \right ] \, dt \,\, ,
\end{eqnarray}
where, in order to obtain the fourth equality of \eqref{calcolo1}, we have used the fact that, since $\beta$ is compactly supported:
\begin{equation}\label{supportocompatto}
0= \int_a^b \,\frac{d}{dt}\left (\frac{\partial L}{\partial \dot{\alpha}} \, \beta \right ) \, dt = \int_a^b \,\left (\frac{d}{dt} \, \frac{\partial L}{\partial \dot{\alpha}}  \right ) \, \beta\, \, dt + \int_a^b \, \frac{\partial L}{\partial \dot{\alpha}} \,\, \dot{\beta} \, \,dt \,\,.
\end{equation}
By way of summary, we conclude from \eqref{calcolo1} that the vanishing of the directional derivative \eqref{derivdir} for all compactly supported variations $\beta$ is equivalent to the Euler-Lagrange equation \eqref{eqeulero1}.
\\

We are now in the right position to extend this setting to the framework of biharmonic maps. In the equivariant context, the bienergy functional $E_2 (\varphi)$ introduced in \eqref{bienergia} takes the following form:
\begin{equation}\label{bienergiaridotta}
    E_2^{\varphi}(\alpha)= \int_a^b \,\, L(t, \alpha, \dot{\alpha}, \ddot{\alpha}) \, dt \,\, ,
\end{equation}
where, as above, the function $\alpha :[a,b] \to \R $ may have to satisfy suitable boundary conditions, and $L(\cdot, \cdot, \cdot, \cdot)$ is a differentiable function depending on the geometry of the problem under consideration. Now, according to Lemma \ref{adaptationofXin}, the condition of biharmonicity for an equivariant map $\varphi_{\alpha}$ with basic bitension field is equivalent to $\alpha $ being a critical point of the reduced bienergy functional \eqref{bienergiaridotta}. Next, we can state our most useful result in this context:

\begin{theorem}\label{reductiontheorem}
A differentiable function $\alpha : [a,b] \to \R $ is a critical point of the reduced bienergy functional \eqref{bienergiaridotta} if and only if it is a solution of the following differential equation:
\begin{equation}\label{eqeulero2}
     \frac{\partial L}{\partial \alpha} - \frac{d}{dt} \, \left ( \frac{\partial L}{\partial \dot{\alpha}} \right ) + \frac{d\,^2}{dt^2} \, \left ( \frac{\partial L}{\partial \ddot{\alpha}} \right )  \, = \, 0 \,\, .
\end{equation}
\end{theorem}

\begin{proof} In \eqref{supportocompatto} we have already noticed that, if $\beta$ is a compactly supported function, then
\begin{equation}\label{supportocompattobis}
  \int_a^b \, \frac{\partial L}{\partial \dot{\alpha}} \,\, \dot{\beta} \, \,dt = - \, \int_a^b \,\left (\frac{d}{dt} \, \frac{\partial L}{\partial \dot{\alpha}}  \right ) \, \beta\, \, dt \,\,.
\end{equation}
In a similar vein, we shall also need the following equality, which holds when both $\beta$ and $\dot{\beta}$ are compactly supported:
\begin{equation}\label{supportocompattotris}
  \int_a^b \, \frac{\partial L}{\partial \ddot{\alpha}} \,\, \ddot{\beta} \, \,dt =  \int_a^b \,\left (\frac{d\,^2}{dt^2} \, \frac{\partial L}{\partial \ddot{\alpha}}  \right ) \, \beta\, \, dt \,\,.
\end{equation}
Now, we proceed to the verification of \eqref{supportocompattotris}: we start from
\begin{equation}\label{supportocompatto4}
0= \int_a^b \,\frac{d}{dt}\left (\frac{\partial L}{\partial \ddot{\alpha}} \, \dot{\beta} \right ) \, dt = \int_a^b \,\left (\frac{d}{dt} \, \frac{\partial L}{\partial \ddot{\alpha}}  \right ) \, \dot{\beta}\, \, dt + \int_a^b \, \frac{\partial L}{\partial \ddot{\alpha}} \,\, \ddot{\beta} \, \,dt \,\,,
\end{equation}
from which we deduce:
\begin{equation}\label{supportocompatto5}
  \int_a^b \, \frac{\partial L}{\partial \ddot{\alpha}} \,\, \ddot{\beta} \, \,dt = - \, \int_a^b \,\left (\frac{d}{dt} \, \frac{\partial L}{\partial \ddot{\alpha}}  \right ) \, \dot{\beta}\, \, dt \,\,.
\end{equation}
Next, we observe that
\begin{eqnarray}\label{supportocompatto6}
0= \int_a^b \frac{d}{dt}\left [ \left (\frac{d}{dt}  \frac{\partial L}{\partial \ddot{\alpha}} \right )  \beta \right ]  dt =  \int_a^b \left (\frac{d\,^2}{dt^2}  \frac{\partial L}{\partial \ddot{\alpha}}  \right )  \beta\, dt + \int_a^b \, \left (\frac{d}{dt}  \frac{\partial L}{\partial \ddot{\alpha}} \right )  \dot{\beta} \, dt. \;\;\;\;
\end{eqnarray}
By using \eqref{supportocompatto6} in \eqref{supportocompatto5} we easily obtain \eqref{supportocompattotris}. We are now in the position to prove the theorem: a function $\alpha$ is a critical point of the functional $E_2^{\varphi}(\alpha)$ in \eqref{bienergiaridotta} if and only if
\begin{equation}\label{derivdir2}
    \frac{d}{dh} \, \left [  E_2^{\varphi}(\alpha \, +\, h\, \beta) \right ] |_{h=0} \, = \, 0
    \end{equation}
for all (differentiable) variations $\beta : [a,b] \to \R$, with $\beta$ and $\dot{\beta}$ compactly supported. Next, we compute the directional derivatives in \eqref{derivdir2} explicitly:
\begin{eqnarray} \label{calcolo2}
\nonumber  \frac{d}{dh} \, \left [  E_2^{\varphi}(\alpha \, +\, h\, \beta) \right ] |_{h=0} &=& \frac{d}{dh}\left ( \int_a^b \,\, L(t, \alpha+h\beta, \dot{\alpha}+h\dot{\beta}, \ddot{\alpha}+h\ddot{\beta})\, dt  \right )\Big {|}_{h=0} \\ \nonumber
   &=& \int_a^b \,\left (\frac{d}{dh}\, L(t, \alpha+h\beta, \dot{\alpha}+h\dot{\beta},\ddot{\alpha}+h\ddot{\beta})|_{h=0} \right )\, dt  \\ \nonumber
   &=& \int_a^b \,\left (  \frac{\partial L}{\partial \alpha} \, \beta +  \frac{\partial L}{\partial \dot{ \alpha}} \, \dot{\beta} + \frac{\partial L}{\partial \ddot{ \alpha}} \, \ddot{\beta}\right ) \, dt \\
   &=& \int_a^b \left [ \left ( \frac{\partial L}{\partial \alpha} - \frac{d}{dt} \,  \frac{\partial L}{\partial \dot{\alpha}} + \frac{d\,^2}{dt^2} \,  \frac{\partial L}{\partial \ddot{\alpha}}\right ) \, \beta \right ] \, dt \,\, ,
\end{eqnarray}
where, in order to obtain the fourth equality of \eqref{calcolo2}, we have used both \eqref{supportocompattobis} and \eqref{supportocompattotris}. Now it follows that the validity of \eqref{derivdir2} for all $\beta$ is equi\-valent to the fact $\alpha$ is a solution of \eqref{eqeulero2}, as required.

\end{proof}

By way of conclusion, an \emph{equivariant} map $\varphi_{\alpha}$ with basic bitension field is \emph{biharmonic} if and only if $\alpha$ is a solution of \eqref{eqeulero2}.

\begin{remark} Our approach require explicitly that the equivariant map $\varphi_{\alpha}$ have basic bitension field. Therefore, it would be interesting to obtain a complete geometric characterization of the situations in which this property is satisfied. In general, this appears to be a rather technical and difficult pro\-blem. However, here we point out that there are important, large families of equivariant maps for which it is immediate to conclude that tension and bitension field are both basic. In particular, that occurs when the Rieman\-nian submersions $\pi_i : E_i \to M_i \,, \,\, i=1,2$ in \eqref{diagrammacommutativo} are determined by isometric actions of  Lie groups, and the equivariant maps $\varphi_{\alpha}$, when restricted to the fibres endowed with the induced metric, are harmonic maps with constant energy density. All the examples in \secref{sezesempi} below are of this type.
\end{remark}
\begin{remark}
For biharmonic curves in a Riemannian manifold, the Euler-Lagrange method was investigated in \cite{CMOP} where, with a different method, an equation of type \eqref{eqeulero2} was derived.
\end{remark}
\begin{remark} In some important geometric applications, such as the study of biharmonic immersions, for instance,
it is necessary to consider a variant of \eqref{bienergiaridotta}, in which the unknown function $\alpha$ is replaced by a curve. More precisely, we have to consider:
\begin{equation}\label{bienergiaridottavettoriale}
    E_2^{\varphi}(\alpha_j)= \int_a^b \,\, L(t, \alpha_j, \dot{\alpha_j}, \ddot{\alpha_j}) \, dt \,\, , \,\quad j= 1 \ldots r \,\, .
\end{equation}
In this case, the argument given in the proof of Theorem \ref{reductiontheorem} applies again and leads us to the conclusion that the critical points of the functional \eqref{bienergiaridottavettoriale} are precisely the solutions of the following \emph{system} of ordinary differential equations:
\begin{equation}\label{eqeulero2vettoriale}
     \frac{\partial L}{\partial \alpha_j} - \frac{d}{dt} \, \left ( \frac{\partial L}{\partial \dot{\alpha_j}} \right ) + \frac{d\,^2}{dt^2} \, \left ( \frac{\partial L}{\partial \ddot{\alpha_j}} \right )  \, = \, 0 \,\, , \,\quad j= 1 \ldots r \,\, .
\end{equation}
\end{remark}

\section{Equivariant biharmonic maps and applications}\label{sezesempi}
In this section we discuss some examples and their applications to problems concerning stability and maximum principle.

\begin{example}\label{esempiotorosfera} {\it Equivariant maps from the flat 2-torus $T^2$ to the 2-sphere $S^2$.}
\end{example}
We write the flat 2-torus as a product manifold:
\begin{equation}\label{toro}
T^2 = \left ( S^1 \times S^1, d\gamma^2+ d\theta^2 \right ) \,\, , \qquad 0 \leq \gamma, \theta \leq 2 \pi \,\,.
\end{equation}
Next, we describe the 2-sphere $S^2$ by means of polar coordinates:
\begin{equation}\label{2sfera}
S^2 = \left ( S^1 \times [0,\pi], \, \sin^2 \alpha \, \, d\beta^2+ d\alpha^2 \right ) \,\, , \qquad 0 \leq \beta\leq 2 \pi \,, \,\, 0 \leq \alpha\leq \pi \,\, .
\end{equation}
We consider equivariant maps $\varphi_{\alpha} : T^2 \to S^2$ of the following form:
\begin{equation}\label{equivdatoroasfera}
    \left ( \gamma, \, \theta \right ) \mapsto \left ( k\, \gamma, \, \alpha (\theta) \right ) \,\, ,
\end{equation}
where $k\in \z$ is a fixed integer, and $\alpha$ is a differentiable, periodic function to be determined (in this case, the period must be equal to $2\pi$). We also note an alternative way to describe an equivariant map of type \eqref{equivdatoroasfera}. More precisely, considering the canonical isometric embedding of $S^2$ into $\R^3$, we can rewrite \eqref{equivdatoroasfera} as follows:
\begin{equation}\label{equivdatoroasferabis}
    \left ( \gamma, \, \theta \right ) \mapsto \left ( (\sin \alpha(\theta) ) \, (\sin k\, \gamma) , \, (\sin \alpha(\theta) ) \, (\cos k\, \gamma), \, (\cos \alpha(\theta) ) \right ) \,\, .
\end{equation}
If one wishes to allow the unknown function $\alpha$ to take values outside the interval $[0, \pi]$, then the form \eqref{equivdatoroasferabis} is to be preferred. In any case, a direct computation shows that the tension field of these equivariant maps can be written as follows:
\begin{equation}\label{tensiontorosfera}
\tau(\varphi_\alpha)= \left [\ddot{\alpha}(\theta) - k^2 \, \sin \alpha(\theta) \,\, \cos \alpha(\theta) \right ] \,\, \frac{\partial}{\partial \alpha} \,\,.
\end{equation}
Therefore, writing $\alpha$ instead of $\alpha (\theta)$, we find that the reduced bienergy \eqref{bienergiaridotta} associated to this family of equivariant maps is given, up to a constant, by:
\begin{eqnarray}\label{bienergiatorosfera}
   E_2^{\varphi}(\alpha) &=&  \int_0^{2 \pi} \,\,  \left [\ddot{\alpha} - k^2 \, \sin \alpha \,\, \cos \alpha \right ] ^2 \,\, d \theta \\ \nonumber
   &=& \int_0^{2 \pi} \,\,  \left [\ddot{\alpha}^2 + \frac{k^4}{4} \, \sin^2 \,(2\alpha) - k^2 \, \ddot{\alpha} \, \sin \, (2 \alpha)\, \right ] \, d \theta \,\, .
\end{eqnarray}
Next, we can compute the condition of biharmonicity by applying directly \eqref{eqeulero2} to \eqref{bienergiatorosfera}. We obtain the following fourth order ordinary differential equation for $ \alpha$:
\begin{equation}\label{biarmoniatorosfera}
    \alpha^{(4)} - \ddot{\alpha} \, \left [ 2 \, k^2 \, \cos (2 \alpha)\right ] + \dot{\alpha}^2 \, \left [ 2 \, k^2 \, \sin (2 \alpha)\right ] + \, \frac{k^4}{2} \, \sin (2\alpha) \, \cos (2 \alpha) \, = \, 0 \,\, .
\end{equation}
At present, as we have already pointed out in the introduction, we do not dispose of general methods to carry out a satisfactory qualitative study of this type of equations. However, in this example we can observe that there are some trivial solutions, namely:
\begin{equation}\label{trivialsolutions}
    \rm {(i)}\,\, \alpha \, \equiv \ell \, \frac {\pi}{2} \,\,, \rm{where} \,\,  \ell =0,\,1,\,2  \, \, ; \qquad \rm {(ii)}\,\, \alpha\, \equiv \frac {\pi}{4} \,\,\, \rm{or}\, \,\,  \alpha \, \equiv \frac {3\, \pi}{4} \,\, .
\end{equation}
The solutions in \eqref{trivialsolutions}(i) correspond to harmonic maps. By contrast, those in \eqref{trivialsolutions}(ii) represent proper biharmonic maps. Although these biharmonic maps were known by other methods (indeed, they can be obtained by projection of $T^2$ onto $S^1$, followed by a biharmonic immersion (see \cite{Montaldo})), our approach leads us to a rather surprising 1-dimensional stability property. More precisely, we can prove the following result:
\begin{theorem}\label{stabilitytheorem} Let $\varphi_{\alpha} : T^2 \to S^2$ be a proper biharmonic map as in \eqref{trivialsolutions}{\rm(ii)}. Then $\alpha$ is a local minimum for the reduced bienergy functional \eqref{bienergiatorosfera}.
\end{theorem}

\begin{proof} We have to compute the second variation of the reduced bienergy functional \eqref{bienergiatorosfera} at the point $\alpha \equiv (\pi \slash 4)$ (which will be denoted by $\alpha^*$). We find:
\begin{eqnarray}\label{secondvariation}
     \nonumber
    \nabla^2 \, E_2^{\varphi}(\alpha^*)\,(V,V) &=&\frac{d\,^2}{dh^2} \, \left [  E_2^{\varphi}(\alpha^* \, +\, h\, V) \right ] |_{h=0} \\ \nonumber
    &\,& \\ \nonumber
    &=& \int_0^{2 \pi} \,\, \frac{d\,^2}{dh^2} \, \left [ (h\, \ddot{V})^2 + \frac{k^4}{4} \, \sin^2 \,(\frac{\pi}{2}+2h\,V) \right .\\ \nonumber
    &\,& \\ \nonumber
    & \,& \left . \qquad \qquad \quad  -h\,\ddot{V}\,k^2 \, \sin \,(\frac{\pi}{2}+2h\,V)\,\right ]\Big |_{h=0} \, \, d \theta \\ \nonumber
    &=& \int_0^{2 \pi} \,\, \frac{d\,^2}{dh^2} \, \left [ h^2 \, \ddot{V}^2 + \frac{k^4}{4} \, \cos^2 \,(2h\,V) \right . \\ \nonumber
    &\,& \\ \nonumber
    &\,& \left . \qquad \qquad \quad - h\,\ddot{V}\,k^2 \, \cos \,(2h\,V)\,\right ] \Big |_{h=0} \, \, d \theta \\
    &=& \int_0^{2 \pi} \, \, \left [ 2 \, \ddot{V}^2 + 2\, V^2 \,k^4 \,\right ] \, d \theta \,\, .
\end{eqnarray}
It follows from the computations in \eqref{secondvariation} that the solution $\alpha^*$ is strictly stable, as required to end the proof. The case $\alpha \equiv (3\, \pi \slash 4)$ is analogous, so we omit it.

\end{proof}

\begin{remark} We point out that the reduced bienergy of the biharmonic maps of Theorem~\ref{stabilitytheorem} is strictly positive. Therefore, they are {\it local}, but clearly {\em not absolute} minima. This fact makes it reasonable to conjecture that Mountain Pass techniques may be used to prove the existence of (unstable) equivariant biharmonic maps of rank 2.
\end{remark}

Now we are going to consider a situation which represents a generalization of our previous Example~\ref{esempiotorosfera}. More precisely, we shall work with warped products of the following type:
\begin{equation}\label{warpedproduct}
    \left (S^m \times (a,b), \, f^2(r)\, g_{S^m} + dr^2 \right ) \,\, ,
\end{equation}
where $\left (S^m , \, g_{S^m}\right )$ denotes the Euclidean unit m-sphere, and $f$ is a smooth, strictly positive function on the interval $(a,b)$. For instance, if $f(r)=r$, then our warped product is (a subset of) the Euclidean space $\R^{m+1}$. The cases $f(r)=\sin r$ and $f(r)=\sinh r$ correspond respectively to $S^{m+1}$ and to hyperbolic space $H^{m+1}$, while $f(r)\equiv c >0$ is a cylinder.

We can introduce the following family of equivariant maps between warped products of type \eqref{warpedproduct}:
\begin{eqnarray}\label{equivarianttrawarped}
 \nonumber
    \varphi_{\alpha}\,: \left (S^m \times (a,b), \, f^2(r)\, g_{S^m} + dr^2 \right ) &\to&
    \left (S^n \times (c,d), \, h^2(\alpha)\, g_{S^n} + d\alpha^2 \right ) \\
    ( \gamma, \, r) \, &\mapsto& \, (\Psi_\lambda(\gamma), \, \alpha(r)) \,\, ,
\end{eqnarray}
where $\Psi_\lambda(\gamma)$ is a so-called {\it eigenmap} of eigenvalue $\lambda$ . That means that
$\Psi_\lambda:\, S^m \to S^n$ is a harmonic map with {\it constant} energy density equal to $(\lambda \slash 2)$. Important examples of eigenmaps are: the identity map of $S^m$ ($\lambda=m$), the k-fold rotation $e^{i\theta} \rightsquigarrow e^{ik\theta}$ of $S^1$ ($\lambda=k^2$); and, also, the Hopf fibrations $S^3 \to S^2$, $S^7 \to S^4$ and $S^{15} \to S^8$, with $\lambda$ equal to 8, 16 and 32 respectively.

Now, the unknown function $\alpha(r)$ in \eqref{equivarianttrawarped} has to be determined in such a way that $\varphi_{\alpha}$ be a biharmonic map. To this purpose, a calculation based on \eqref{tensionfield} shows that the tension field of an equivariant map of the form \eqref{equivarianttrawarped} is given by the following expression:
\begin{equation}\label{tensionfieldtrawarped}
    \tau (\varphi_{\alpha})= \left [\ddot{\alpha}(r) + m\, \frac{\dot{f}(r)}{f(r)}\, \dot{\alpha}(r)- \lambda \, \frac{h(\alpha)\,\dot {h}(\alpha)}{f^2(r)} \right ] \,\, \frac{\partial}{\partial \alpha} \,\,.
\end{equation}
Therefore, in this case the reduced bienergy \eqref{bienergiaridotta} is given, up to an irrelevant constant factor, by the following expression:
\begin{equation}\label{reducedbienergiatrawarped}
    E_2^{\varphi}(\alpha) =  \int_a^b \,\,  \left [\ddot{\alpha}(r) + m\, \frac{\dot{f}(r)}{f(r)}\, \dot{\alpha}(r)- \lambda \, \frac{h(\alpha)\,\dot {h}(\alpha)}{f^2(r)}\right ] ^2 \,\,f^m(r)\, dr \,\, .
\end{equation}
Now, a direct application of \eqref{eqeulero2} to \eqref{reducedbienergiatrawarped} yields the biharmonicity equation for $\alpha$. But, since this equation is in general rather long and complicated, we shall write it down and analyze in a few particular instances only. We start with:
\begin{example}\label{esempiowarpedeuclidean} {\it Equivariant maps between Euclidean spaces} (compare with \cite{Montaldobis}).
\end{example}
This situation corresponds to the case that $f(r)=r$ and $h(\alpha)=\alpha$ in \eqref{reducedbienergiatrawarped}. If we also perform the following change of variable:
\begin{equation}\label{cambiovariabile euclideanspaces}
    r=e^t \, , \quad t \in \R \, , \qquad \beta(t)= \alpha(e^t) \, \, ,
\end{equation}
then the reduced bienergy functional \eqref{reducedbienergiatrawarped} takes the following form:
\begin{equation}\label{reducedbienergiabeta}
    E_2^{\varphi}(\beta) =  \int_{\R} \,\,  \left [\ddot{\beta}+ (m-1)\, \dot{\beta}\, - \lambda \,\beta\right ] ^2 \,\,e^{(m-3)t} \, \, dt \,\, .
\end{equation}
By way of example, we assume that $\Psi_\lambda: S^3 \to S^2$ is the Hopf map, so that $m=3$ and $\lambda =8$. A simple computation shows that the biharmonicity equation \eqref{eqeulero2} becomes:
\begin{equation}\label{biarmonicitatraeuclideihopf}
    \beta^{(4)} -20 \, \ddot{\beta} + 64 \, \beta =0 \,\, ,
\end{equation}
which admits the following explicit representation of solutions:
\begin{equation}\label{biarmonicheesplicitehopfeuclideo}
    \beta(t)= c_1 \, e^{4t}+c_2 \, e^{-4t}+c_3 \, e^{2t}+c_4 \, e^{-2t} \,\, .
\end{equation}
We observe that, when either $c_1$ or $c_4$ is different from zero, the map which is associated to \eqref{biarmonicheesplicitehopfeuclideo} is \emph{proper biharmonic}: this fact is actually part of a general principle (Almansi's property), which states that multiplication of a harmonic function on $\R^{m+1}$ by the factor $r^2$ generates a proper biharmonic function (see \cite{Montaldo}).
\begin{example}\label{esempiocilindrotoeuclidean} {\it Equivariant maps from cylinders to Euclidean spaces (violation of the maximum principle).}
\end{example}
Now we illustrate the case of equivariant maps $\varphi_{\alpha}: \,S^m \times \R \to \R^{n+1}$. This instance corresponds to the choice $f(r)\equiv 1, \, r \in \R $ and $h(\alpha)=\alpha$ in \eqref{reducedbienergiatrawarped}, so that the reduced bienergy functional \eqref{reducedbienergiatrawarped} becomes:
\begin{equation}\label{reducedbienergiadacilindroaeuclideo}
    E_2^{\varphi}(\alpha) =  \int_{\R} \,\,  \left [\ddot{\alpha}- \lambda \, \alpha) \right ] ^2 \,\, dr \,\, .
\end{equation}
Next, applying \eqref{eqeulero2} to \eqref{reducedbienergiadacilindroaeuclideo}, we find that the biharmonicity equation for this case is:
\begin{equation}\label{biarmonicitadacilidroaeuclideo}
    \alpha^{(4)} -2 \, \lambda \, \ddot{\alpha} + \lambda^2 \, \alpha =0 \,\,
\end{equation}
(note the similarity between equation \eqref{biarmonicitadacilidroaeuclideo} and a time independent version of \eqref{trave}).

Equation \eqref{biarmonicitadacilidroaeuclideo} admits the following explicit representation of solutions:
\begin{equation}\label{biarmonicheesplicitedacilindroaeuclideo}
    \alpha(r)= c_1 \, e^{\sqrt{\lambda}\, r}+c_2 \, e^{-\sqrt{\lambda}\, r}+c_3 \, r \, e^{\sqrt{\lambda}\, r}+c_4 \, r \, e^{-\sqrt{\lambda}\, r} \,\, .
\end{equation}
We observe that, when either $c_3$ or $c_4$ is different from zero, the map associated to \eqref{biarmonicheesplicitedacilindroaeuclideo} is a proper biharmonic map of full rank $(n+1)$ (of course, provided that we are using an eigenmap of rank $n$).

In \cite{Sampson}, Sampson proved the following maximum principle for harmonic maps (to state this result, we keep the notation of \cite{Sampson}):
\begin{theorem}\label{maxprinctheorem} Assume that $f:M \to Y$ is a harmonic map, with $q=f(p)$. Let $S$ be a $C^2$ hypersurface in $Y$ passing through $q$, at which point we assume that the second fundamental form is definite. If $f$ is not a constant mapping, then no neighbourhood of $p$ is mapped entirely on the concave side of $S$.
\end{theorem}

In order to construct our counterexample, we consider an equivariant biharmonic map arising from the following solution of \eqref{biarmonicheesplicitedacilindroaeuclideo} ($\lambda >0$ , $\Psi_\lambda$ of full rank):
\begin{equation}\label{esempiopermaxprinciple}
    \alpha(r)= (\sqrt {\lambda} \, r ) \, \sinh (\sqrt {\lambda} \, r ) \, + \, e^{(\sqrt {\lambda} \, r )} \,\, , \qquad r \in \R \,\, .
\end{equation}
It is easy to verify that the solution $\alpha(r)$ in \eqref{esempiopermaxprinciple} admits a \emph{strictly positive absolute minimum point}, say $r=r_0$. Therefore, the image of our biharmonic map is entirely contained in the concave side of $S=\partial B_{r_0}(O)$ . This fact shows that proper biharmonic maps do not satisfy the Sampson maximum principle of Theorem \ref{maxprinctheorem}.

\end{document}